\documentclass[12pt]{article}
\title{Braid Groups and Right Angled Artin Groups}
\author{Frank Connolly and Margaret Doig}
\date{\today}
\usepackage{amsmath,amsthm,amscd,amssymb}
\chardef\bslash=`\\ 
\theoremstyle{definition}

\newtheorem{thm}{Theorem}[section]

\newtheorem{con}[thm]{Construction}
\newtheorem{definition}[thm]{Definition}
\newtheorem{lemma}[thm]{Lemma}

\newtheorem{prop}[thm]{Proposition}

\newtheorem{guess}[thm]{Conjecture}
\newcommand{\p}{{\partial}}
\newcommand{\cR}{{\cal R}}
\newcommand{\cT}{{\cal T}}

\swapnumbers

\begin{document}
\maketitle

\newcommand{\U}{\cal U \rm _n^{top}(X)}
\newcommand{\Ul}{\cal U\rm_{n-1}^{top}(X)}


\begin{abstract} In this article we prove a special case of a
conjecture of A. Abrams and R. Ghrist about fundamental groups of
certain aspherical spaces. Specifically, we show that the $n-$point
braid group of a {\em linear tree} is a right angled Artin group for each n.
\end{abstract}


\section{Introduction. Statement of Results.}

A  \emph{graph} is a  connected one dimensional compact
polyhedron.  Ghrist and   Abrams (\cite{abramsthesis},
\cite{abramsgeom},\cite{ghrist}) have recently called attention to
the $n-$point unordered configuration space of a graph $X$,
denoted here as $\U$. This is the space of $n-$element subsets of
$X$ (see Definition \ref{unx}, (\ref{confspace})). It is an
aspherical space with the  homotopy type of a  finite polyhedron,
for each $n$ and $X$ (see \cite{ag}). Its fundamental group is the
\emph{$n-$string braid group of $X$,} denoted $B_n(X,c)$, if $c$ is
a base point of $\U$. The group $B_n(X,c)$ is therefore torsion free; it
can have arbitrarily high finite cohomological dimension. Abrams and
Ghrist (\cite{abramsgeom}, \cite{ghrist}) have put forward the
following striking conjecture:

\guess \label{conj} If $X$ is a planar graph, then the $n-$string
braid group $B_n(X,c)$ is a right angled Artin group for each $n$.

A \emph{right angled Artin group} is a group having a presentation
in which the only relations are commutators between generators. It
is known (\cite{abramsgeom}) that the planar condition cannot be
removed from the above conjecture.

The purpose of this paper is to prove:

\thm  \label{raathm}For each $n$,  Conjecture \ref{conj} is true  if
$X$ is a linear tree.

A \emph{linear tree} is a contractible graph $X$ containing an
interval $I$ (a homeomorphic copy of $[0,1]$) such that each node
of $X$ is in $I$.

(Recall that the \emph{nodes} in a graph $X$ are the points of
degree $\geq 3$. The \emph{degree} of any point  $x\in X$ is the
number of points in $Link_X(x)$).

Theorem \ref{raathm} is a direct consequence of our main theorem
(\ref{mainthm}) below. Before stating it we give a more careful
definition of the objects just mentioned.

\begin{definition} \label{unx}  A subset $c$ with $n$ elements in a
space $X$ is called an $n-$point configuration in $X$. For any
$n\geq 0$ and any topological space $X$, the unordered $n-$point
configuration space of $X$ is:
\begin{equation}\label{confspace} \U = \{ c \subset X| \quad |c|=n\} .
\end{equation}
If $(X,d)$ is a metric space, then $\U$ is topologized by using
the Hausdorff metric on closed sets of $X$.  Explicitly, then:
\begin{equation}
d(c,c') = max\{d(x,c'), d(y,c) |\quad x\in c, \;y\in c'\}\quad
\forall  c,c'\in \U.
\end{equation}
Equivalently, one can  give $ \U$ the quotient topology of the map
\[(X^n-\Delta) \overset{\pi}{\longrightarrow} \U : \quad \pi(x_1,
x_2, \dots x_n)= \{x_1,x_2,\dots x_n\}\] where $\Delta=\{(x_1,
x_2, \dots x_n)| \quad x_i= x_j  \text{ for some } i\neq j\}$.

Note \;$ \cal{U}\rm^{top}_0(X)= \ast, $\; and  $\cal U \rm
_1^{top}(X) = X$.

The \emph{$n-$string braid group} of a space $X$ is:
\begin{equation}\label{braid}
B_n(X, c) = \pi_1( \U, c),
\end{equation}
where $c\in \U$ is a base point.
\end{definition}

$B_n$ is a functor from the category of topological spaces $X$
with fixed base configurations $c\in \U$ and isotopy classes of
injective continuous maps that preserve base configurations.

But in order to finesse the many changes of base configuration
required, we  employ the following artifice.  Let $(X,A)$ be a
pair where $A$ is a nonempty {\em simply connected} subspace of
$X$. The \emph{fundamental group} of $(X,A)$ denoted $\pi_1(X,A)$
is the  set of homotopy classes  of maps $(I, \p I)
\overset{\sigma}{\longrightarrow} (X,A)$

The multiplication in this group is:
\[
[\sigma][\tau] = [\sigma\cdot \rho\cdot \tau]\text{ where $\rho$
is any path in $A$ from $\sigma(1)$ to $\tau(0)$}
\]
This group is functorial in such $(X,A)$ and isomorphic to
$\pi_1(X,x_0)$ if $x_0\in A$. In particular, if $I $ denotes an
interval (a subspace homeomorphic to $[0,1]$) in $X$, then $\cal
U\rm_n^{top}(I)$ is a  contractible subset of $\U$. We define the
\emph{n-point braid group of (X,I)} by:
\begin{equation}\label{braid XI}
B_n(X;I) = \pi_1( \U, \cal U\rm_n^{top}(I)).
\end{equation}

\begin{definition} \label{i_p}(The Endpoint Inclusion Map). Let $X$ be
a graph. Let $p$ be an endpoint of $X$ (that is, $degree(p) =1$).
For each $n\geq 1$, we define a map
\[ B_{n-1}(X,c)\overset{\iota_p}{\longrightarrow} B_n(X,c') \]
as follows. Choose an isotopy of $1_X$, say $\{r_t:X\to X
\;\;|\;\; 0\leq t\leq 1\}$, which is stationary outside a small
neighborhood of $p$ and which satisfies $r_1(X)\subset X-\{p\}.$
Then $r_1$ induces a map
\[\Ul\overset{\i_p}{\longrightarrow}\U:\quad i_p(d) = \{p\}\cup
r_1(d) \quad \forall \;d\in \Ul.
\]
In turn, $i_p$ induces a map of fundamental groups, denoted:
\[B_{n-1}(X,c)\overset{\iota_p}{\longrightarrow} B_n(X,c'), \text{ if } c'=\iota_p (c)
\]
(or $\iota_{p\;X}$ when such explicitness is needed).

Now if $I$ is an interval in $X$ containing $p$ and if the isotopy
fixes $X-I$, then $i_p$ induces
\[B_{n-1}(X;I)\overset{\iota_p}{\longrightarrow} B_n(X;I)\]
which is independent of the isotopy chosen. Write
$\iota_p^k:B_{n-k}(X;I)\to B_n(X;I)$ for the $k$-fold iteration of
this map. Abrams (\cite{abramsgeom}, Lemma 3.4, p.24) shows that
$\iota_p$ is injective if $X$ is any graph.
\end{definition}
(Note: Abrams proves this for the {\em pure} braid groups. These
have finite index in our braid groups $B_n(X;I)$, which  are
torsion free. This implies that $\iota_p$ is injective).

Our main theorem will say that the groups $B_n(X;I)$, for
$n=0,1,2,\dots$, admit right angled Artin presentations which are
all related by the maps $\iota_p$ above. A \emph{right angled
Artin presentation} of a group $G$, denoted $\langle \beta,\cal R
\rangle $, consists of a subset $\beta$ of $G$, and a subset $\cR$
of $F(\beta)$, the free group  on $\beta$, such that $\cR$
consists of elements of the form $(xyx^{-1}y^{-1})$, where $x,y\in
\beta$, and the following sequence is exact
\[
1\to N\to F(\beta)\overset{j}{\longrightarrow}G\to 1,
\]
where $j$ denotes the natural  homomorphism, and $N$ denotes the
smallest normal  subgroup
containing $\cR$.

Here is the main theorem.

\thm \label{mainthm}Let $X$ be a tree. Let $p$  be an endpoint of
$X$. Let $I\subset X$ be an interval containing $p$ and every node
of $X$.  Then for each integer $n\geq0$ there is a right angled
Artin presentation, $\langle \beta(n), \cR(n)\rangle$ for
$B_n(X;I)$ such that
\begin{equation}
\iota_{p}(\beta(n-1))\subset \beta (n) \quad \text{and }
\iota_{p\ast}(\cR(n-1))\subset \cR(n)\quad \forall n\geq 1.
\label{nat}
\end{equation}
Here $\iota_{p\ast}:F(\beta(n-1))\to F(\beta(n))$ is the
homomorphism induced by the function $ \iota_p:\beta(n-1)\to \beta(n)$.

\

It is easy to see that Theorem \ref{raathm} follows from Theorem
\ref{mainthm} because every interval in a tree $X$ lies in a
bigger interval $I$ containing an endpoint of $X$.

\

Here is an outline of the rest of this paper. In
Section~\ref{sec2}, we study the case of a \emph{star} (a tree
with one node). If $X$ is a star and $I$ is an interval whose
endpoints are endpoints of $X$, we show that $B_n(X;I)$ is a free
group admitting a basis $\beta(n)$, for each $n$, such that
$\iota_p(\beta(n-1))\subset \beta(n)$ if $p$ is \emph{either}
endpoint.  In Section \ref{sec3}, we prove a theorem computing the
braid group $B_n(X;I)$ if $X$ is the one point union of two graphs
along a common endpoint. In section \ref{sec4}, we use the
previous results to prove the main theorem, \ref{mainthm}.

We wish to thank Aaron Abrams here for suggesting this problem to us, in 2003. His comments
about it were encouraging and helpful.

\section{Braid Groups of Stars.}\label{sec2}

A \emph{star} $S$ is a tree with no more than one node. The
node is denoted $v$. If $S$ is a tree with no nodes (an interval), then
we assume further that a fixed interior point $v$ of $S$ is given.
It turns out that, for any star $S$ and any integer $n\geq 0$, the
$n-$point configuration space $\cal U\rm_n^{top}(S)$ contains a
compact one dimensional polyhedron $D_n(S)$ which is a deformation
retract of $S$.  This is constructed in \cite{mdoig}. We review
the construction  here. We  then use $D_n(S)$ to prove the
following:

\prop \label{starprop} Let $S$ be a star. Let $p$ and $q$ be two
distinct endpoints of $S$, and let $I$ be the interval $[p, q]$.
Then for each $n\geq 0$,  $B_n(S;I)$ is a free group. Moreover,
there is  a basis $\beta(n)$ for $B_n(S;I)$ such that, if $n\geq
1$,
\[ \iota_p(\beta(n-1))\subset \beta (n) \text{ and}
\quad\iota_q(\beta(n-1))\subset \beta (n). \]

Note: The rank of this free group $B_n(S;I)$ is
\[
1+(k-1)\binom{n+k-2}{k-1} - \binom{n-k-1}{k-1}
\]
where $k$ is the number of endpoints of $S$. This is proved by
Doig in \cite{mdoig}, but an equivalent formula for the
corresponding \emph{pure}, braid group of $S$  appears earlier in
Ghrist \cite{ghrist}, Prop. 4.1.

\

Before beginning the proof of \ref{starprop}, we now choose a
fixed metric $d$ on the star $S$. To this end, fix the integer
$n\geq0$. The metric $d$ will be a constant multiple of the
standard simplicial metric $\rho$ on $S$. The star $S$ has a {\em
canonical} simplicial structure in which the vertices are the
points of degree different from $2$.  For the corresponding
simplicial metric $\rho$, $[p, q]$ is isometric to $[0,2]$ and $v$
is the midpoint of $[p, q]$. Define the metric $d$ by:
\begin{equation}\label{metric}
d(x,y) = C\rho(x,y) \quad \forall \, x,y\in S
\end{equation}
where $C$ is a fixed constant such that $C\geq max(1, n-1) $.

\

\con\label{DnS}(of $D_n(S)$; compare \cite{mdoig}).

Let $S$ be a star. Let \;$ c\in\cal U\rm_n^{top}(S)$ be any
$n-$point configuration in the star $S$. For each endpoint $p$ of
$S$, we define
\begin{equation}\label{arm}
A_p(c)= c\cap [p,v].
\end{equation}
We call this an \emph{arm} of $c$ \rm if $A_p(c)\neq \emptyset$.

We say $c$ is \emph{regular} if it satisfies the following rules:
\begin{enumerate}\label{regular}
\item  For all $x,y\in c $ with $x\neq y, \; d(x,y)\geq 1.$ Also
$d(x,y) =1$ if $x$ and $y $ \;lie in a single arm of $c$ and
$[x,y]\cap c=\{x, y\}$. \item If $v\notin c$, and if $A_p(c)$ is
an arm of $c$, then there is a point $x\in c$ such that $d(x,
A_p(c)) =1$.
\end{enumerate}
In English: successive points in a single arm of $c$ are one unit
apart; if $v\notin c$, the innermost point of each arm of $c$ has
distance one from some other arm, and has distance at least one
from every other arm.

Therefore $c$ has  at least one  arm $A_p(c)$ such that
$d(v,A_p(c))\leq \frac{1}{2} $. There is {\em at most} one arm of
$c$ that satisfies: $0<d(v,A_p(c))<\frac{1}{2}$  (the
\emph{governing arm} in the language of \cite{mdoig}). When
$A_p(c)$ is the unique governing arm of $c$, every other arm,
$A_q(c)$, satisfies: $d(v, A_q(c))= 1-d(v, A_p(c))$.

The subspace  $D_n(S) $ of $\cal U\rm_n^{top}(S)$ can now be
defined:
\begin{equation}\label{dns}
D_n(S) := \{c\in \cal U\rm_n^{top}(S)| \; c \text{ is regular}\}.
\end{equation}
Note $D_0(S) = \ast =\cal U\rm_0^{top}(S)$, and  $D_1(S) =\{v\}$.

$D_n(S)$ has the structure of a 1-dimensional compact polyhedron,
as we now explain. Its vertices come in two types. The Type I
vertices are those configurations $c\in D_n(S)$ such that $v\in c
$.  The Type II vertices  are those $c\in D_n(S)$ such that
$d(v,c)=\frac{1}{2}$. For Type  II vertices, note that
$d(v,A_p(c))=\frac{1}{2}$ for each arm $A_p(c)$.

Each 1-simplex of $D_n(S)$ has a single Type I vertex $c$ and a
single Type II vertex $c'$. We  denote this 1-simplex $[c,c']$ or
$[c',c]$. For each Type  II vertex $c'$, there is a single
1-simplex $[c', c]$ for each endpoint $p$ of $S$ such that
$A_p(c')\neq \emptyset$. The other endpoint $c$ of this 1-simplex
$[c',c]$ is defined as the unique Type I vertex such that:
\begin{equation}\label{type one vertex} |A_p(c)| = |A_p(c')|,\textrm{ and
}|A_q(c)| = |A_q(c')|+1 \;\; \text{for any \emph{other} endpoint }
q. \end{equation} The  points of the 1-simplex $[c', c]$ are $c'$,
$c$, and those $e\in D_n(S)$ such that
\begin{equation} |A_q(e)| = |A_q(c')| \;\text{ for all endpoints
$q$ of $S$, and } 0<d(v, A_p(e))<\frac{1}{2}. \end{equation} The
rule $e\mapsto d(v,A_p(e))$ gives a homeomorphism from $[c',c]$
onto the interval $[0,\frac{1}{2}]$.

Each point $e\in D_n(S)$, other than a vertex, belongs to a unique
1-simplex $[c', c]$. The Type II vertex $c'$ is specified by the
requirement that $|A_q(e)|=|A_q(c')| $ for all endpoints $q$. The
endpoint $p$, and also (by (\ref{type one vertex}))  the Type~I
vertex $c$, \;is specified by the requirement that \;
$0<d(v,A_p(e))<\frac{1}{2}$.

This completes the construction of the 1-dimensional polyhedron
$D_n(S)$.

\

Doig proves in \cite{mdoig} that $D_n(S)$ is a  strong deformation
retract of $\cal U\rm_n^{top}(S)$.

Incidentally, if we change the metric $d$ on $S$ to $d'$, by
changing the constant $C$ to $C'$, then $D_n(S, d)$ is isometric
to $D_n(S, d')$. But the two are not identical.

Before beginning the proof of Proposition \ref{starprop}, we need
to relate $\iota_p$ to the deformation retract $D_n(S)$.

If $p$ and $q$ are two endpoints of $S$ and $I=[p, q]$, then
$D_n(I)\subset  D_n(S)$ and $D_n(I)$ is an interval.

Analogous to the map $i_p$ of Definition \ref{i_p} is a simplicial
inclusion map $ \tilde{i}_p:D_{n-1}(S)\to D_n(S)$, defined by:
\[\tilde{i}_p(c)=\{x\}\cup c \quad \quad\forall c\in D_{n-1}(S),
\] where $x$ is the unique point of $[p, v]$ such that $\{x\}\cup c \in D_n(S)$.
(We have chosen the constant $C$ above, so that $d(p,c)\geq 1.$ This
ensures that there {\em is} such a point $x$). It is elementary to see
that
$\tilde{i}_p(D_n(I))\subset D_n(S)$, and that the following diagram 
commutes up to homotopy:
\[ \begin{CD}
(D_{n-1}(S), D_{n-1}(I))@>\tilde{i}_p>> (D_n(S), D_n(I))\\
@V\text{incl.}VV  @V\text{incl.}VV\\
(\cal U\rm^{top}_{n-1}(S), \cal U\rm^{top}_{n-1}(I))@>{i}_p>> (\cal
U\rm^{top}_{n}(S), \cal U\rm^{top}_{n}(I))
\end{CD}
\]
Therefore we can identify  the group $B_n(S;I)$  with
$\pi_1(D_n(S), D_n(I))$, and we can identify the homomorphism
$\iota_p: B_{n-1}(S; I)\to B_n(S;I)$ with the map
$(\tilde{i}_p)_*: \pi_1(D_{n-1}(S), D_{n-1}(I))\to \pi_1(D_n(S),
D_n(I))$.

\begin{proof} (of Proposition \ref{starprop}):

Suppose a maximal tree $\cT$ of $D_n(S)$ is chosen, containing
$D_n(I)$. Because $dim(D_n(S))= 1$, $B_n(S;I)$ is a free group,
and a basis for $B_n(S;I)$ is given by those 1-simplices of
$D_n(S)$ which  are not in $\cT$. Therefore, it is enough to
exhibit, for each $n$, a maximal tree $\cT(n)$ for $D_n(S)$ such
that
\begin{equation}\label{T} \cT(n)\supset D_n(I) \quad \text{
and } \quad \cT(n-1)= \tilde{i}_p^{-1}(\cT(n))=
\tilde{i}_q^{-1}(\cT(n)).
\end{equation}
We do this  now.

Let $c^{(n)}$ be  the unique point of $D_n(S)$ such that
$c^{(n)}\subset [p, v]$. The configuration $c^{(n)}$ is a Type  I  vertex.
For each vertex $c$ of $D_n(S)$ except $c^{(n)}$, we are going to
construct a \emph{successor vertex} $s(c)$  such that $[c, s(c)]$ is a
1-simplex and $s^k(c) = c^{(n)}$ for some $k>0$. Then $ \cT(n)$
will consist of the union of these $[c, s(c)].$

Number the endpoints of $S: \; p_1, p_2,\dots p_m$. Ensure that
$p=p_1$ and $q=p_2$. Define $s(c)$ as follows. If $c$ is a Type I
vertex, then $s(c)$ is the unique Type II vertex satisfying:
\begin{equation}\label{scI}
|A_{p_1}(s(c))| = |A_{p_1}(c)|; \quad\qquad |A_{p_j}(s(c))|=
|A_{p_j}(c)|-1 \text{ if } j\neq 1
\end{equation}
If $c$ is a Type II vertex, let $r=r(c)$ be the biggest index for
which $A_{p_{r(c)}}\neq \emptyset$. Note $r\geq 2 $ (by
\ref{regular}.2.). Define $s(c)$ as the unique Type I vertex
satisfying:
\begin{equation}\label{scII}
|A_{p_r}(s(c))| = |A_{p_r}(c)|;\qquad |A_{p_j}(s(c))| =
|A_{p_j}(c)|+1 \text{ if } j\neq r.
\end{equation}
For any vertex $c$  of $D_{n-1}(S)$ except $c^{(n-1)}$, we have, by
(\ref{scI}) and
(\ref{scII})
\begin{equation}\label{scIII} \tilde{i}_{p_j}(s(c)) =
s(\tilde{i}_{p_j}(c)),  \text { for }
j\leq 2.
\end{equation}
By (\ref{type one vertex}), $c$ and $s(c)$ span a 1-simplex for each
vertex $c\neq c^{(n)}$.
It is also clear from   (\ref{scI})  and (\ref{scII}) that, for each
vertex $c$, there is an
integer
$k\geq 0$ such that $s^k(c) = c^{(n)}$.

Therefore we define:
\begin{equation}\label{Tn}
\cT(n) = \cup\{[c, s(c)]\; \;|\; c \text{ is a vertex of } D_n(S)
\text{ other than }
c^{(n)}\}.
\end{equation}
$\cT (n)$ is a tree containing every vertex of $D_n(S);$ therefore it
is maximal.

By (\ref{scIII}), $\; \tilde{i}_{p_j}(\cT(n-1))\subset \cT(n)$, for
$j=1$ and $2.$  But since $\cT(n-1)$ is a \emph{maximal} tree  in
$D_{n-1}(S)$, and $\cT(n)$ contains no cycles, this
  implies:
\[
\cT(n-1) = \tilde{i}_p^{-1}(\cT(n))= \tilde{i}_q^{-1}(\cT(n)).
\]
Finally we must show that $D_n(I)\subset \cT(n)$. The key is to
note that each Type II vertex $c$ of $D_n(I)$ belongs to exactly
two 1-simplices. One of these is $[c, s(c)]$. The other is  $[c,
d]$, where $d$ is that Type I vertex of $D_n(I)$ such that
$c=s(d)$. It follows that $D_n(I)\subset \cT(n)$. This completes
the proof of Proposition \ref{starprop}.
\end{proof}

\section{The Endpoint Union of Two Graphs.}\label{sec3}

This section is devoted to computing the $n$-point braid group of the
union of two graphs which intersect at a single endpoint of each.

Our goal is Proposition~\ref{exact} below.

We let $X$ be a graph of the form $X=S \cup T$ where $S$ and $T$
are graphs.  Assume that $\{q\} = S \cap T$, a single point, and
that $q$ is an endpoint of $S$ and of $T$.  Let
$j:S \longrightarrow X $ and $ j': T\longrightarrow
X,$ be  inclusion maps. These induce maps of braid groups with the same
names.   Let
$I$ and
$J$ be intervals in
$S$ and $T$ respectively so that $\{q\}=I \cap J$. Form the free
product $B_n(S;I) \ast B_n(T;J)$.  Let $N$ be the smallest normal
subgroup of this product containing each of the  commutator subgroups
$[\iota_q^k B_{n-k}(S;I),\iota_q^{n-k} B_k(T;J)],~1<k<n$.  Then we have:

\begin{prop}\label{exact}
With the hypotheses above, the following sequence is exact: \[1
\rightarrow N \rightarrow B_n(S;I) \ast B_n(T;J) \stackrel{j \ast
j'}{\longrightarrow} B_n(X;I \cup J) \rightarrow 1.\] 
\end{prop}

We will need:

\begin{lemma} \label{exactlem}
Let $A_0 \subseteq A_1 \subseteq \dots \subseteq A_n$ and $B_0
\subseteq B_1 \subseteq \dots \subseteq B_n$ be two increasing
sequences of groups.  Suppose that the following diagram is a
pushout diagram: 
\[ G 
\]
 \[
\nearrow \qquad \qquad \qquad \uparrow \qquad \qquad \qquad \nwarrow 
\]
\[\qquad A_n\times B_0 \qquad \qquad \qquad A_{n-1}\times B_1\qquad
\qquad\qquad A_{n-2}\times B_2\quad \dots
\]
\[
  \nwarrow   i\times 1    \qquad    \nearrow 1\times i'\qquad 
\nwarrow i\times 1 \quad  \nearrow  1\times  i'
\]
\[
  A_{n-1}\times B_0  \qquad \qquad \qquad A_{n-2}\times
B_1\quad \dots 
\]
 
(i.e., $G$ is the direct limit of the diagram obtained by
deleting
$G$ and the maps to $G$). Here $i$ and $i'$ denote inclusions.
  Then there is an exact sequence:
 \[1
\rightarrow N \rightarrow A_n \ast B_n \stackrel{j \ast
j'}{\longrightarrow} G \rightarrow 1\] where $N$ is the smallest
normal subgroup of $A_n*B_n $ containing $[A_{n-k},B_{k}]$ for each
$k=0,1,2,
\dots, n$. Here $j$ and  $j'$ are restrictions of the limit maps
$A_n\times B_0\to G$, and $ A_0 \times B_n\to G$ to $A_n$ and $B_n$
respectively.
\end{lemma}

\begin{proof}
If $n=1$ this is clear.  Working inductively, we let $G'$ be the
limit of the diagram obtained by omitting $A_0$ and $ B_n$.  This
gives us $G'=A_n \ast B_{n-1} / N'$ where  $N'$ is the smallest
normal subgroup containing all $[A_{n-k},B_k]$ for $k=0,1,2,
\dots, n-1$. It also  gives a new diagram  whose limit is
obviously $A_n \ast B_n /N$ where $N$ is the smallest normal
subgroup containing $N'$ and $[A_0,B_n]$.
\end{proof}

\begin{proof} ({\em of Proposition~\ref{exact}}):
For each $k=0,1,\dots,n$, let
 \[U(k)=\{c \in \U|\quad|c \cap S|
\geq k, |c \cap T| \geq n-k\},\] 
and let $j_k:U(k)\to  \U$ denote the inclusion map. Since $|S \cap T|=1$,
we see that 
$ U(j) \cap U(k) = \emptyset$ if $|k-j|>1$ and that $\U=\bigcup_{k=0}^n
U(k)$.  Each $c\in U(k)$ can be written uniquely in the form:
\[c=c_S(k) \cup c_T(n-k)\] where:\newline (i) $c_S(k) \subset
S$;\quad (ii) $c_T(n-k) \subset T$;\quad (iii) $|c_S(k)|=k;\quad
(iv) |c_T(n-k)|=n-k$. Necessarily, $c_S(k) \cap
c_T(n-k)=\emptyset$.

One sees easily that there is a homotopy equivalence
 \[U(k)
\stackrel{h_k}{\longrightarrow} \cal U\rm_k^{top}(S) \times \cal
U\rm_{n-k}^{top}(T): \qquad h_k(c)= \left(c_S(k),c_T(n-k)
\right).\] It sends $U(k) \cap \cal U\rm_n^{top}(I \cup J)$ to
$\cal U\rm_k^{top}(I) \times \cal U\rm_{n-k}^{top}(J)$. Similarly
\[U(k)\cap U(k-1) = \{c \in \U~|~q \in c, |c \cap S|=k, |c
\cap T|=n-k+1\}.
\] Therefore there is a homotopy equivalence for each $k\geq 1$:
\[U(k) \cap U(k-1) \stackrel{j'_k}{\longrightarrow}\cal U\rm_{k-1}^{top}(S) \times
\cal U\rm_{n-k}^{top}(T):\quad j_k'(c)=  \left(c_S(k-1),c_T(n-k)
\right).
\]

Passing to fundamental groups of these spaces, we obtain the
pushout diagram below by using the version of Van Kampen's Theorem
in Hatcher (\cite{hatcher} 1.20, p.43):\newpage
\[ 
\qquad  B_n(X;I\cup J)\]
\[  
j_0 \nearrow  \qquad \qquad  \qquad  j_1\uparrow   \qquad \qquad
\qquad  \nwarrow j_2
\]
\[B_0(S;I)\times B_n(T;J)  \quad    \quad B_1(S;I)\times B_{n-1}(T;J) 
\quad   \quad B_2(S;I)\times B_{n-2}(T;J)\dots\]
 \[ \nwarrow i\times 1  \qquad 1\times i' \nearrow   \qquad \qquad 
\nwarrow i\times 1 \qquad 1\times  i' \nearrow   \]
\[  B_0(S;I)\times B_{n-1}(T;J)  \qquad \qquad   B_1(S;I)\times
B_{n-2}(T;J)\dots  \]

\

  Now, by
Lemma~\ref{exactlem}, where $j=j_n$, $j'=j_0$, we get the exact
sequence \[1 \rightarrow N \rightarrow B_n(S;I) \ast B_n(T;J)
\overset{j*j'}\rightarrow B_n(X;I \cup J) \rightarrow 1.\]
\end{proof}

\section{Proof of the Main Theorem.}\label{sec4}
\begin{proof}
(of the Main Theorem \ref{mainthm}): If the tree $X$ has less than
two nodes, the proof is already clear from
Proposition~\ref{starprop}. Therefore we  assume $X$ has at least
two nodes. We write $X$ as: $X= S\cup T$, where $S$ is a star with one
node, $v$, $T$ is a tree,  $S\cap T  =\{q\}$,  where $q$ ia an endpoint
of both $S$ and $T$, and 
$p$ is an endpoint of $X $  lying in $S$.
 We can
moreover assume $I=J\cup K$, where $J= [p,q]$ and $K$ is an interval  in
$T$,  containing $q$ and each node of $T$.

 From Proposition
\ref{starprop} we have a basis $\beta_S(n)$ for the free group
$B_n(S;J)$ for each $n$, such that
\begin{equation}\label{betan}
\iota_p(\beta_S(n-1))\subset \beta_S(n);\quad
\iota_q(\beta_S(n-1))\subset \beta_S(n).
\end{equation}
By induction on the number of nodes of $X$, we have a right angled
Artin presentation    $\langle \beta_T(n); \cR_T(n)\rangle $ of
$B_n(T;K)$ for  each  $n= 0, 1, 2, \dots$ so that
\begin{equation}
\label{artinT} \iota_q(\beta_T(n-1)) \subset \beta_T(n),\qquad
\iota_{q*}(\cR_T(n-1))\subset \cR_T(n), \quad \forall  n\geq 1.
\end{equation}
By Proposition \ref{exact} we have an exact sequence \[1\to
N\rightarrow B_n(S;J)*B_n(T;K) \overset{j_S*j_T} {\longrightarrow}
B_n(X;I)\to 1 \] where $j_S$ and $j_T$ are induced by the
inclusions $S\subset X$ and $ T\subset X$  and $N$ is the smallest
normal subgroup of the free product containing each of the sets:
\[ [\iota_q^k\beta_S(n-k), \,\iota_q^{n-k}\beta_T(k)] . \]
Set \begin{align} \beta_X(n) &= j_S(\beta_S(n))\cup j_T(\beta_T(n))\label{artpres1}\\
\cR_X(n) &= j_{T*}(\cR_T(n))\cup [j_{S*}\iota^{k}_{q*}\beta_S(n-k),
\; j_{T*}\iota^{n-k}_{q*}\beta_T(k)]\label{artpres2}
\end{align} where $j_{S*}:F(\beta_S(n))\to F(\beta_X(n))$ and
$j_{T*}:F(\beta_T(n))\to F(\beta_X(n))$ are induced by $j_S$ and
$j_T$.

Clearly $\langle \beta_X(n); \cR_X(n)\rangle$ is an Artin
presentation of $B_n(X;I)$ for each $n$. To complete the
argument, we  must prove (\ref{nat}). First note that
\begin{align} \iota_{p X}\circ j_T &= j_T\circ \iota_{q T}: B_{n-1}(T;K)\to B_n(X;
I)\label{comm1}\\
\iota_{p X}\circ j_S  &= j_S\circ \iota_{p S}: B_{n-1}(S;J)\to B_n(X;
I)\label{comm2}\\
\iota_{p S}\circ \iota_{q S} &= \iota_{q S}\circ \iota_{p S}: B_{n-2}(S;J)\to
B_n(S;J)\label{comm3}
\end{align} because the corresponding diagrams of spaces commute up to homotopy.
It follows that
\begin{equation}\label{ibx} \iota_p(\beta_X(n-1))\subset \beta_X(n)
\end{equation}
by (\ref{artpres1}). Also, by (\ref{artinT}), (\ref{artpres2}) and
(\ref{ibx}) we have
\begin{equation} 
(\iota_{p \, X})_*j_{T*}
(\cR_T(n-1)) = j_{T*}(\iota_{q\, T})_*(\cR_T(n-1))\subset
j_{T*}(\cR_T(n))\subset \cR_X(n).\label{last}
\end{equation}
Finally by ( \ref{comm1}), (\ref{comm2}) and (\ref{comm3}) we have:
\begin{multline*} \iota_{p\;X}[j_{S\;*}(\iota^k_{q\;S})_*\beta_S(n-1-k),\quad
j_{T\;*}(\iota_{q\;T}^{n-1-k})_*\beta_T(k)]\\
\subset [j_{S\;*}(\iota^k_{q\;S})_*\beta_S(n-k),
j_{T\;*}(\iota_{q\;T}^{n-k})_* \beta_T(k)]
\end{multline*}
which, with  (\ref{ibx}) and (\ref{last}) implies \[(\iota_{p
X})_*\cR_X(n-1) \subset \cR_X(n). \] This completes the proof of
Theorem \ref{mainthm}.
\end{proof}

\end{document}